\documentclass[11pt,a4paper,oneside]{article}
\usepackage{epsfig}
\usepackage{amsfonts}
\usepackage[latin1]{inputenc}
\usepackage{amsmath}
\usepackage{amssymb}

\newtheorem{teo}{\underline{Theorem}}[section]
\newtheorem{cor}{\underline{Corollary}}[section]
\newtheorem{prop}{\underline{Proposition}}[section]
\newtheorem{lem}{\underline{Lemma}}[section]
\newcommand{\R}{\mathbb{R}}
\newcommand{\pr}{\mathbb{P}}

\newcommand{\be}{\begin{equation}}
\newcommand{\ee}{\end{equation}}

\title{Polynomial rate convergence to an invariant measure for the continuum time limit of the Minority Game}
\author{Matteo Ortisi\thanks{UniCredit Markets\&Investment Banking, via Broletto 16, 20121 Milano, Italy.
matteo.ortisi@gmail.com}
}

\begin{document}
\maketitle

\begin{abstract}
In this paper we show that the continuum time version of the  Minority Game satisfies the criteria for the application of a theorem on the existence of an invariant measure. We consider the special case of a game with ``sufficiently" asymmetric initial condition where the number of possible choices for each individual is $S=2$ and $\Gamma<+\infty$. An upper bound for the asymptotic behavior, as the number of agents grows to infinity, of the waiting time for reaching the stationary state is then obtained.
\end{abstract}

\section{Introduction} 
The Minority Game (MG) \cite{challet} is a simple model based on Arthur's ``El Farol" bar problem \cite{arthur} which describes the behavior of a group of competing heterogeneous agents subject to the economic law of supply and demand. An application of the MG is, for example, the microscopic modeling of financial markets \cite{challet_fin1, challet_fin, marsililibro, laureti}.

In this paper the attention is focused on the continuum time version of the MG (see perhaps \cite{cavagna, garrahan}), where the number of possible choices for each individual is $S=2$. In particular we are interested in studying its long time behavior, in terms of existence of an invariant measure for its dynamical variables and convergence of the dynamical variables distribution to it.

Let us fix some notations, coherent with those used in \cite{marsili1}.
Consider the MG with $N$ agents. Its dynamics is defined in terms of dynamical variables $U_{s,i}(t)$ in discrete time $t=0,1,\ldots$; these are scores corresponding to each of the possible agents choices $s=+1,-1$. Each agent takes a decision $s_i(t)$ with
$$
Prob\{s_i(t)=s\}=\frac{e^{\Gamma_i U_{s,i}(t)}}{\sum_{s'}e^{\Gamma_iU_{s',i}(t)}}
$$
where $\Gamma_i>0$. The original MG corresponds to $\Gamma_i=\infty$ and was generalized to $\Gamma_i=\Gamma<\infty$ \cite{cavagna}; here we consider the latter case.

The public information variable $\mu(t)$ is given to all agents; it belongs to the set of integers $\{1,\ldots,P\}$ and can either be the binary encoding of last $M$ winning choices or drawn at random from a uniform distribution; here we consider the latter case.

The action $a_{s_i(t),i}^{\mu(t)}$ of each agent depends on its choice $s_i(t)$ and on $\mu(t)$. The coefficients $a_{s,i}^{\mu}$, called strategies, are uniform random variables taking values $\pm 1$ ($Prob\{a_{s,i}^{\mu}=\pm 1 \}=1/2$) independent on $i, s$ and $\mu$.

Let us introduce the following random variables (to ease the notation the choices $+1$ and $-1$ are shorted with $+$ and $-$)
$$\xi_i^{\mu}=\frac{a_{+,i}^{\mu}-a_{-,i}^{\mu}}{2}, \quad \Theta^{\mu}=\sum_{i=1}^N\frac{a_{+,i}^{\mu}+a_{-,i}^{\mu}}{2}$$
and their averages
$$
\overline{\xi_i\Theta}=\frac{1}{P}\sum_{\mu=1}^P\xi_i^{\mu}\Theta^{\mu}, \quad \overline{\xi_i\xi_j}=\frac{1}{P}\sum_{\mu=1}^P\xi_i^{\mu}\xi_j^{\mu}.
$$
The only relevant quantity in the dynamics is the difference between the scores of the two strategies:
$$
y_i(t)=\Gamma\frac{U_{+,i}(\tau)-U_{-,i}(\tau)}{2},
$$
where $\tau=\frac{t}{\Gamma}$.

Let $(\Omega,\mathcal{F},\pr)$ be the probability space respect to which all our random variables are defined, $y=(y_i)_{1\leq i\leq N}$, $\Theta=(\Theta^{\mu})_{1\leq\mu\leq P}$ and $\xi=(\xi_i)_{1\leq i\leq N}$.

As shown in \cite{marsili1}, if $P/N=\alpha\in\R_+$,  $S=2$ and $\Gamma_i=\Gamma>0$ for all $i$, the dynamics of the continuum time limit of the MG is given by the following $N$-dimensional stochastic differential equation
\begin{equation}
\label{stoch_eq1}
dy_i(t)=\left(-\overline{\xi_i\Theta}-\sum_{j=1}^N\overline{\xi_i\xi_j}\tanh (y_j)\right)dt + A_i(y,N,\Gamma,\xi)dW(t), \quad i=1,\ldots,N
\end{equation}
where
\begin{itemize}
\item[]$W(t)$ is an $N$-dimensional Wiener process,
\item[]$A_i$ is the $i$-th row of the $N\times N$ matrix $A=(A_{ij})$ such that $$(AA')_{ij}(y,N,\Gamma,\xi)=\frac{\Gamma\sigma_{N,\Gamma}^2(y)}{\alpha N}\overline{\xi_i\xi_j}.$$
\end{itemize}
If $\alpha>\alpha_c$, where $\alpha_c=0.3374\ldots$ marks the transition point from a symmetric $(\alpha<\alpha_c)$ to an asymmetric phase $(\alpha>\alpha_c)$ characterized, respectively, by no predictability and predictability (in the asymmetric phase the choices $+1$ and $-1$ do not appear with equal probabilities for a given $\mu(t)$), then the function  $\sigma_{N,\Gamma}^2:\R^N\to\R_+$ is continuous and 
$$
\lim_{N\to\infty}\sup_{y\in\R^N}\frac{\sigma_{N,\Gamma}^2(y)}{N}\leq 1,
$$ 
(see \cite{cavagna, marsili1}).

In \cite{marsili1} the authors derive the full stationary distribution of $y$. Here, by proving that equation (\ref{stoch_eq1}) satisfies the criteria for the application of Veretennikov's Theorem (see \cite{veretennikov}), we show for the continuum time version of the MG with $S=2$, ``sufficiently" asymmetric initial condition, $\Gamma<+\infty$ and $\alpha>\alpha_c$, that, for $N$ sufficiently large, the distribution $\nu(t)$ of $y$ converges almost surely with polynomial rate to an invariant measure $\nu_{inv}$ and that the waiting time for reaching the stationary state is at most  $O(\vert y(0)\vert^2)$. For finite asymmetric initial condition this means that the waiting time is at most $O(N)$ or, since $\frac{P}{N}=\alpha>0$, $O(\frac{P}{\alpha})$.  

The work is organized as follows. In section \ref{preliminary} some preliminary computations useful to characterize the asymptotic behavior of the drift coefficient of equation (\ref{stoch_eq1}) are performed. In section \ref{invariant} Veretennikov's Theorem is applied to equation (\ref{stoch_eq1}) with both finite and maximally ($\vert y(0)\vert\to\infty$) asymmetric initial condition and an upper bound for limiting behavior as $N$ grows to infinity of the waiting time for reaching the stationary state is obtained. In section \ref{conclusions} some conclusions are drawn,  while Appendix contains the Veretennikov's Theorem.

\section{Preliminaries}
\label{preliminary}
Let us start by proving some preliminary results showing that, despite the fact that the coefficients present in the drift term of equation (\ref{stoch_eq1}) are random variables assuming both negative and positive values, as $N$ grows to infinity the behavior of equation (\ref{stoch_eq1}) becomes dissipative and hence suitable to the application of a stability theorem, like the Veretennikov's one. 
For the sake of simplicity in notation, let us assume $P=N$, i.e. $\alpha=1$. The results obtained still hold for any $\alpha>0$.
From now on, when defining a probability event we shall omit to explicit the dependence of the random variables on the $\omega$s $\in \Omega$.

\begin{lem}
\label{lemma1}
For every $i=1,2,\ldots$
\begin{eqnarray*}
\label{prob_inequality}
&&\pr\left\{\lim_{N\to\infty}\overline{\xi_i\Theta}=0 \wedge \lim_{N\to\infty}\overline{\xi_i^2}=a \wedge \lim_{N\to\infty}\sum_{j=1, j\neq i}^N\overline{\xi_i\xi_j}=0  \right\}=\nonumber\\
&&  \pr\left\{ \lim_{N\to\infty}\overline{\xi_i^2}=a\right\} ,\nonumber\\
\end{eqnarray*}
where $a$ is any constant.
\end{lem}
\begin{flushleft}
Proof
\end{flushleft}
Obviously, for every $i$
\begin{eqnarray*}
&&\pr\left\{\lim_{N\to\infty}\overline{\xi_i\Theta}=0 \wedge \lim_{N\to\infty}\overline{\xi_i^2}=a \wedge \lim_{N\to\infty}\sum_{j=1, j\neq i}^N\overline{\xi_i\xi_j}=0  \right\}\leq\\
&&  \pr\left\{ \lim_{N\to\infty}\overline{\xi_i^2}=a\right\}.
\end{eqnarray*}

It remains hence to show that also the inverse inequality holds.

For every $i$,
\begin{eqnarray}
\label{prob_inequality1}
\!\!\!\!\!\!&&\pr\left\{\lim_{N\to\infty}\overline{\xi_i\Theta}=0 \wedge \lim_{N\to\infty}\overline{\xi_i^2}=a \wedge \lim_{N\to\infty}\sum_{j=1, j\neq i}^N\overline{\xi_i\xi_j}=0  \right\}=\nonumber\\
\!\!\!\!\!\!&&\pr\left\{\lim_{N\to\infty}\left[\frac{1}{N}\sum_{\mu=1}^N\left[\frac{(a_{+,i}^{\mu})^2-(a_{-,i}^{\mu})^2}{4}+\xi_i^{\mu}\sum_{j=1, j\neq i}^N\frac{a_{+,j}^{\mu}+a_{-,j}^{\mu}}{2}\right]\right]=0 \wedge \lim_{N\to\infty}\overline{\xi_i^2}=a\right.\nonumber\\
\!\!\!\!\!\!&& \left.\wedge  \lim_{N\to\infty}\sum_{j=1, j\neq i}^N\frac{1}{N}\sum_{\mu=1}^N \xi_i^{\mu}\xi_j^{\mu}=0  \right\}.
\end{eqnarray}
Since $\{\xi_i^{\mu}\}_{\mu,i}$ is a family of independent (on both $\mu$ and $i$) random variables identically distributed  over the set $\{-1,0,1\}$  and since $(a_{+,i}^{\mu})^2-(a_{-,i}^{\mu})^2=0$, the right hand side term of (\ref{prob_inequality1}) is greater than or equal to
\begin{eqnarray*}
\label{prob_inequality2}
\pr\left\{\lim_{N\to\infty}\frac{1}{N}\sum_{\mu=1}^N\!\sum_{j=1, j\neq i}^N\!\!\!\frac{a_{+,j}^{\mu}+a_{-,j}^{\mu}}{2}=0 \wedge \lim_{N\to\infty}\overline{\xi_i^2}=a
\wedge  \lim_{N\to\infty}\!\!\!\sum_{j=1, j\neq i}^N\!\!\frac{1}{N}\sum_{\mu=1}^N \xi_j^{\mu}=0  \right\}=\nonumber\\
\pr\left\{ \lim_{N\to\infty}\overline{\xi_i^2}=a\right\} \pr\left\{\lim_{N\to\infty}\frac{1}{N}\sum_{\mu=1}^N\!\sum_{j=1, j\neq i}^N\!\!\!\!\!\frac{a_{+,j}^{\mu}+a_{-,j}^{\mu}}{2}=0 \wedge \!\!  \lim_{N\to\infty}\!\!\!\sum_{j=1, j\neq i}^N\!\!\frac{1}{N}\sum_{\mu=1}^N \xi_j^{\mu}=0  \right\}.\nonumber\\
\end{eqnarray*}
Let $\eta_j=\sum_{\mu=1}^N\xi_j^{\mu}$ and $\zeta_j=\sum_{\mu=1}^N\frac{a_{+,j}^{\mu}+a_{-,j}^{\mu}}{2}$ ($\eta_j$ and $\zeta_j$ both depend on $N$, but since for every $N$, $\mathbb{E}[\eta_j]=\mathbb{E}[\zeta_j]=0$ we omit to explicitly write the $N$ dependence); $\left\{\eta_j+\zeta_j=\sum_{\mu=1}^N a_{+,j}^{\mu}\right\}_j$ is a family of independent identically distributed random variables with mean $\mathbb{E}[\eta_j+\zeta_j]=0$.

Since by the Law of Large Numbers (LLN)
$$
\pr\left\{ \lim_{N\to\infty}\frac{1}{N}\sum_{j=1, j\neq i}^N(\eta_j+\zeta_j)=0 \right\}=1 \quad \pr\left\{ \lim_{N\to\infty}\frac{1}{N}\sum_{j=1, j\neq i}^N\eta_j=0 \right\}=1,
$$
it follows that
$$
\pr\left\{ \lim_{N\to\infty}\frac{1}{N}\sum_{j=1, j\neq i}^N\zeta_j=0 \right\}=1
$$
and the thesis follows.

\hspace{12cm}$\Box$

\begin{lem}
\label{lemma2}
For every $i=1,2,\ldots$
$$
\pr \left\{\lim_{N\to\infty}\overline{\xi_i^2}=\frac{1}{2}\right\}=1.
$$
\end{lem}    
\begin{flushleft}
Proof
\end{flushleft}
Since $\mathbb{E}[(\xi_i^{\mu})^2]=\frac{1}{2}$, it is a consequence of the LLN.

\hspace{12cm}$\Box$

\begin{lem}
\label{useful_prop}
For every $i=1,2,\ldots$
$$
\pr\left\{\lim_{N\to\infty}\overline{\xi_i\Theta}=0 \wedge \lim_{N\to\infty}\overline{\xi_i^2}=\frac{1}{2} \wedge \lim_{N\to\infty}\sum_{j=1, j\neq i}^N\overline{\xi_i\xi_j}=0  \right\}=1.
$$
\end{lem}
\begin{flushleft}
Proof
\end{flushleft}
It is an immediate consequence of Lemma \ref{lemma1} and Lemma \ref{lemma2}.

\hspace{12cm}$\Box$

\begin{lem}
\label{prop_2.2}
Let
\begin{eqnarray*}
B_i&=&\left\{ \lim_{N\to\infty}\overline{\xi_i^2}=\frac{1}{2} \bigwedge_{j=1,2,\ldots, j\neq i} \lim_{N\to\infty}\overline{\xi_i\xi_j}=0 \right\}\\
i&=&1,2,\ldots
\end{eqnarray*}
Then
$$
\pr\left( \bigcap_{i=1,2,\ldots}B_i \right)=1.
$$
\end{lem}
\begin{flushleft}
Proof
\end{flushleft}
For every $i=1,2,\ldots$
\begin{eqnarray}
\label{result1}
\pr\left(B_i\right)&\geq& \pr\left\{ \lim_{N\to\infty}\overline{\xi_i^2}=\frac{1}{2} \bigwedge_{j=1,2,\ldots, j\neq i} \lim_{N\to\infty}\frac{1}{N}\sum_{\mu =1}^N \xi_j^{\mu}=0 \right\}\nonumber\\
&=&\pr\left\{ \lim_{N\to\infty} \overline{\xi_i^2}=\frac{1}{2}\right\}\prod_{j=1,2,\ldots, j\neq i}\pr\left\{\lim_{N\to\infty} \frac{1}{N}\sum_{\mu=1}^N\xi_{j}^{\mu}=0 \right\}=1,
\end{eqnarray}
where last equality is due to Lemma \ref{lemma2} and the LLN.

Let
$$
D_1=B_1, \quad D_i=B_i\cap B_{i-1}, \quad i=2,3,\ldots
$$
Obviously $\pr\left( D_1 \right)=\pr\left( B_1 \right)=1$. Under the inductive hypothesis
$\pr\left( D_{i-1} \right)=1$; by (\ref{result1})
$$
\pr\left(D_i\right)=\pr\left(D_{i-1}\cap B_i\right)=\pr\left(B_i\vert D_{i-1}\right)\pr\left(D_{i-1}\right)=1.
$$
Since $\left\{D_i\right\}_{i\in\mathbb{N}}$ is a decreasing sequence 
$$
\pr\left( \bigcap_{i=1,2,\ldots}B_i \right)=\pr\left( \bigcap_{i=1,2,\ldots}D_i \right)=\lim_{i\to\infty}\pr\left( D_i \right)=1.
$$

\hspace{12cm}$\Box$


Let us define the following events:
\begin{eqnarray*}
A_1\!\!\!&=&\!\!\!\left\{ \lim_{N\to+\infty}\lim_{x_1\to\pm\infty}\frac{1}{N}\left[\overline{\xi_1\Theta}x_1+\overline{\xi_1^2}\tanh(x_1)x_1+\sum_{j=2}^N\overline{\xi_1\xi_j}\tanh(x_j)x_1\right]=\frac{1}{2} \right\},\\
A_i\!\!\!&=&\!\!\!\left\{ \lim_{N\to+\infty}\lim_{\vert x\vert\to\infty}\frac{1}{N}\left[\overline{\xi_i\Theta}x_i+\overline{\xi_i^2}\tanh(x_i)x_i+\sum_{j=1, j\neq i}^N\overline{\xi_i\xi_j}\tanh(x_j)x_i\right]\in\left[0,\frac{1}{2}\right] \right\}\\
i&=&2,3\ldots
\end{eqnarray*}
and 
$$
E_1=A_1,\quad
E_i=A_i\cap E_{i-1}, \quad i=2,3,\ldots
$$

\begin{lem}
\label{lemma_A_i}
For every $i=1,2,\ldots$,
$$
\pr\left(A_i\right)=1.
$$
\end{lem}
\begin{flushleft}
Proof
\end{flushleft}

Since
$$
\lim_{N\to+\infty}\lim_{x_1\to\pm\infty}\frac{\tanh(x_1)x_1}{N}=1\quad {\rm and}\quad \lim_{x_1\to\pm\infty, x_j\to\pm\infty}\frac{\tanh(x_1)x_1}{\tanh(x_j)x_1}= \pm 1
$$
by Lemma \ref{useful_prop}, $\pr\left(A_1\right)=1$.
About to the $A_i$, i=2,3,\ldots, if  $x_i\to\pm\infty$, obviously $\pr\left(A_i\right)=1$.
If $x_i\neq\pm\infty$, since $\tanh$ is a bounded odd function, by Lemma \ref{useful_prop}, $\pr\left(A_i\right)=1$.

\hspace{12cm}$\Box$

\begin{lem}
\label{lemma_E_i}
For every $i=1,2,\ldots$,
$$
\pr\left(E_i\right)=1.
$$
\end{lem}
\begin{flushleft}
Proof
\end{flushleft}
The proof is the same of Lemma \ref{prop_2.2} with $E_i$ and $A_i$ instead of $D_i$ and $B_i$.

\hspace{12cm}$\Box$

\section{Convergence to the Invariant Measure}
\label{invariant}
In this section we show that, for $N$ sufficiently large, the distribution of the random variable $y$, the score differences vector whose dynamics is described by the stochastic differential equation (\ref{stoch_eq1}), admits an invariant measure and we study the limiting behavior of the  waiting time for reaching the stationary state.

For this purpose we perform an opportune rescaling of the variable $y$ and show that the dynamics of the rescaled variable $z$ satisfies the criteria for the application of Veretennikov's Theorem (Theorem \ref{teo_vere} in Appendix or see \cite{veretennikov}). As a second step we extend the thesis of Veretennikov's Theorem to the original random variable $y$.

The Veretennikov's Theorem gives, under quite general regularity assumptions, a condition that suffices the existence of an invariant measure and the convergence to it for the distribution of a random variable satisfying a stochastic differential equation. 
The criteria for the application of such a theorem to a stochastic differential equation are based on the evaluation of the drift term on the initial condition that must be such that $\vert y(0)\vert\neq 0$; this means that we must avoid a game where all the agents have symmetric initial condition (i.e. $y(0)=0$). We focus our attention on two cases: a game where there exists at least one agent having maximally asymmetric initial condition ($\vert y_i(0)\vert \to\infty$) and no matter about the other agents (they can or not have symmetric initial condition) and a game where $\vert y(0)\vert<\infty$ and the number of agents with asymmetric initial condition ($y_i(0)\neq 0$) is $O(N)$. The former situation corresponds to a game where there is at least one agent playing always the same strategy since the beginning (a so called producer), while the latter one corresponds to a game where  the number of agents perceiving a strategy more successful that the other one since the begging is $O(N)$.

Let $b_i^N(\xi,\Theta,x):(\{-1,0,1\}^{P\times N},\R^P,\R^N)\to\R$ defined as
\begin{eqnarray*}
b_i^N(\xi,\Theta,x)&=&\overline{\xi_i\Theta}+\sum_{j=1}^N\overline{\xi_i\xi_j}\tanh(x_j),\\
i&=&1,2,\ldots,N
\end{eqnarray*}
and
$$
b^N(\xi,\Theta,x)=\left(b_i^N(\xi,\Theta,x)\right)_{1\leq i\leq N}.
$$
Equation (\ref{stoch_eq1}) can be written in the form
\begin{equation}
\label{eq_y}
dy(t)=-b^N(\xi,\Theta,y)dt+A(y,N,\Gamma,\xi)dW(t).
\end{equation}
$b^N(\xi,\Theta,y)$ is obviously a Borel-measurable locally bounded function.

Let $c$ be a bounded constant (that may depend on $N$) greater than $1$ and $z=cy$; under this rescaling of the variable $y$ equation (\ref{eq_y}) becomes
\begin{equation}
\label{eq_z}
dz(t)=-cb^N\left(\xi,\Theta,\frac{z}{c}\right)dt+cA\left(\frac{z}{c},N,\Gamma,\xi\right)dW(t).
\end{equation}
Before going to study the limiting behavior of equation (\ref{eq_z}) as $N$ grows to infinity, let us prove the nondegeneracy of the the diffusion matrix $A$, a condition for the application of Veretennikov's Theorem.

\begin{prop}
\label{prop_3.2}
There exists $\hat{N}>0$ such that for every $N>\hat{N}$ the diffusion matrix $A$ is almost surely nondegenerate.
\end{prop}
\begin{flushleft}
Proof
\end{flushleft}
Nondegenerancy is equivalent to the following condition
$$
\inf_y\inf_{\xi}\inf_{\vert x\vert=1}xAA'(y,N,\Gamma,\xi)x'>0 
$$
where $\Gamma$ is fixed and $x\in\R^N$.
 
Since $\sigma^2_{N,\Gamma}>0$ it is sufficient to show that
$$
\pr\left\{
\lim_{N\to\infty}
\left(
\begin{array}{cccc}
\overline{\xi_1^2} & \overline{\xi_1\xi_2} & \cdots &\overline{\xi_1\xi_N}\\
\overline{\xi_2\xi_1} & \overline{\xi_2^2} & \cdots &\overline{\xi_2\xi_N}\\
\vdots & \vdots & \vdots & \vdots\\
\overline{\xi_N\xi_1} & \overline{\xi_N\xi_2} & \cdots &\overline{\xi_N^2}
\end{array}
\right)
=
a\mathbf{I}
\right\}=1
$$
where $a$ is a positive constant and $\mathbf{I}$ is the identity matrix.

This follows from Lemma \ref{prop_2.2} and hence we obtain the thesis.

\hspace{12cm}$\Box$

\subsection{Minority Game with maximally asymmetric initial condition}
We consider a game where there exists at least one producer, that is an agent $i$ such that $\vert y_i(0)\vert\to\infty$.
We show that, under an appropriate choice of the constant $c$, system (\ref{eq_z}) satisfies the criteria for the application of Veretennikov's Theorem and then we extend the results of Theorem \ref{teo_vere} to the original system (\ref{eq_y}).
To ease the notation let $y(0)=x$.
\begin{prop}
\label{main_proposition}
$$
\pr\left\{\lim_{N\to\infty}\lim_{\vert x\vert\to\infty}\frac{\left\langle b^N(\xi,\Theta,x),x\right\rangle}{N}\geq \frac{1}{2}\right\}=1.
$$
\end{prop}
\begin{flushleft}
Proof
\end{flushleft}
\begin{eqnarray*}
&&\pr\left\{\lim_{N\to\infty}\lim_{\vert x\vert\to+\infty}\frac{\left\langle b^N(\xi,\Theta,x),x\right\rangle}{N}\geq \frac{1}{2}\right\}=\nonumber\\
&&\pr\left\{\lim_{N\to +\infty}\lim_{\vert x\vert \to+\infty}\frac{1}{N}\sum_{i=1}^{N}\left[\overline{\xi_i\Theta}x_i+\overline{\xi_i^2}\tanh(x_i)x_i+\!\!\!\!\!\!\sum_{j=1, j\neq i}^N\overline{\xi_i\xi_j}\tanh(x_j)x_i\right]\geq\frac{1}{2}\right\}\nonumber=\nonumber\\
&&\pr\left\{\sum_{i=1}^{\infty}\lim_{N\to +\infty}\lim_{\vert x\vert \to+\infty}\frac{1}{N}\left[\overline{\xi_i\Theta}x_i+\overline{\xi_i^2}\tanh(x_i)x_i+\!\!\!\!\!\!\sum_{j=1, j\neq i}^N\overline{\xi_i\xi_j}\tanh(x_j)x_i\right]\geq\frac{1}{2}\right\}\nonumber\geq\nonumber\\
&&\pr\left(\bigcap_{i=1,2,\ldots}A_i\right)=\pr\left(\bigcap_{i=1,2,\ldots}E_i\right).
\end{eqnarray*}
Since $\left\{E_i\right\}_{i\in\mathbb{N}}$ is a decreasing sequence,
$$
\pr\left(\bigcap_{i=1,2,\ldots}E_i\right)=\lim_{i\to+\infty}\pr\left(E_i\right).
$$
By applying Lemma \ref{lemma_E_i} the thesis follows.

\hspace{12cm}$\Box$

\begin{cor}
\label{cor_3.1}
There exist $\tilde{N}>0$ and $M_0>0$ such that for every $N>\tilde{N}$
$$
\left\langle b^N(\xi,\Theta,x),\frac{x}{\vert x\vert}\right\rangle> \left(\frac{N}{2}-1\right)\frac{1}{\vert x\vert} \quad \vert x\vert\geq M_0 \quad a.s.
$$
\end{cor}
\begin{flushleft}
Proof
\end{flushleft}
It is an immediate consequence of Proposition \ref{main_proposition}.

\hspace{12cm}$\Box$

\begin{prop}
\label{criterium_prop}
Equation (\ref{eq_z}) with $c>1+\frac{2}{\frac{N}{2}-1}$ satisfies a.s. condition (\ref{cond_vere}) with $r=c\left(\frac{N}{2}-1\right)>\frac{N}{2}+1$.
\end{prop}
\begin{flushleft}
Proof
\end{flushleft}
It is an immediate consequence of Corollary \ref{cor_3.1} and Proposition \ref{prop_3.2}.

\hspace{12cm}$\Box$

In Proposition \ref{main} Veretennikov's Theorem is applied to equation (\ref{eq_z}), while the existence of an invariant measure for $y$ and the rate of convergence to it is derived in Corollary \ref{cor1}.

\begin{prop}
\label{main}
Let $\nu(t)$ be the distribution of $z(t)=(z_i(t))_{1\leq i\leq N}$ satisfying equation (\ref{eq_z}) with initial condition $z(0)$ and $c>1+\frac{2}{\frac{N}{2}-1}$. For $N$ sufficiently large, there exists a.s. an invariant measure $\nu_{inv}$ for $z$ such that
$$
\vert \nu(t)-\nu_{inv}\vert\leq m(1+\vert z(0)\vert^l)(1+t)^{-(k+1)},
$$
where $m$ is a positive constant, $0<k<c\left(\frac{N}{2}-1\right)-\frac{N}{2}-1$, $2k+2<l<2c\left(\frac{N}{2}-1\right)-N$ and $\vert \nu(t)-\nu_{inv}\vert$ is the total variation distance between $\nu(t)$ and $\nu_{inv}$, i.e. $\vert \nu(t)-\nu_{inv}\vert:=\sup_{A\in\mathcal{B}_{\R^d}}\vert \nu(t)(A)-\nu_{inv}(A)\vert$.
\end{prop}
\begin{flushleft}
Proof
\end{flushleft}
It is a consequence of Proposition \ref{criterium_prop} and Theorem \ref{teo_vere}.

\hspace{12cm}$\Box$

\begin{cor}
\label{cor1}
Let $\nu'(t)$ be the distribution of $y(t)=(y_i(t))_{1\leq i\leq N}$ satisfying equation (\ref{eq_y}) with initial condition $y(0)$ and let $c=1+\frac{2}{\frac{N}{2}-1}+\frac{2}{\left(\frac{N}{2}-1\right)^2}$. For $N$ sufficiently large, there exists a.s. an invariant measure $\nu'_{inv}$ for $y$ such that
\begin{equation}
\label{inequality}
\vert \nu'(t)-\nu'_{inv}\vert\leq m'(1+\vert cy(0)\vert^l)(1+t)^{-(k+1)},
\end{equation}
where $m'$ is a positive constant, $0<k<\frac{2}{\frac{N}{2}-1}$ and $2k+2<l<2+\frac{4}{\frac{N}{2}-1}$.
\end{cor}
\begin{flushleft}
Proof
\end{flushleft}
Since the limit as $N$ grows to infinity of $c$ is 1 and $z=cy$, for $N$ sufficiently large $\vert \nu'(t)-\nu'_{inv}\vert=m'\vert \nu(t)-\nu_{inv}\vert$ and the thesis follows from Proposition \ref{main}.

\hspace{12cm}$\Box$

\subsection{Minority Game with finite asymmetric initial condition}

Now we move to a game where $\vert y(0)\vert<\infty$ and the number of agents with asymmetric initial condition ($y_i(0)\neq 0$) is $O(N)$.

Let $y(0)=x$ and $N'=\#\left\{y_i : y_i(0)\neq 0\right\}$ and suppose that $N'=O(N)$, i.e. $\lim_{N\to\infty}\frac{N'}{N}=\gamma$, $0<\gamma\leq 1$; moreover let $\beta_i=\tanh(x_i)x_i$ and $\beta=\min\{\beta_i : \beta_i\neq 0\}$. To ease the notation, let us suppose that $x_1$ is such that $\beta=\tanh(x_1)x_1$.
\begin{prop}
\label{main_proposition_no_frozen}
$$
\pr\left\{\lim_{N\to\infty}\frac{\left\langle b^N(\xi,\Theta,x),x\right\rangle}{N}\geq \frac{1}{2}\beta\gamma\right\}=1.
$$
\end{prop}
\begin{flushleft}
Proof
\end{flushleft}
\begin{eqnarray*}
&&\pr\left\{\lim_{N\to\infty}\frac{\left\langle b^N(\xi,\Theta,x),x\right\rangle}{N}\geq \frac{1}{2}\beta\gamma\right\}=\\
&&\pr\left\{\lim_{N'\to\infty}\frac{\gamma}{N'} \sum_{i=1}^{N'}\left[\overline{\xi_i\Theta}x_i+\overline{\xi_i^2}\tanh(x_i)x_i+\!\!\!\!\!\!\sum_{j=1, j\neq i}^{N'}\overline{\xi_i\xi_j}\tanh(x_j)x_i\right]\geq \frac{1}{2}\beta\gamma\right\}\geq\\
&&\pr\left\{\lim_{N'\to\infty}\left[\overline{\xi_1\Theta}x_1+\overline{\xi_1^2}\tanh(x_1)x_1+\!\!\!\!\!\!\sum_{j=1, j\neq 1}^{N'}\overline{\xi_1\xi_j}\tanh(x_j)x_1\right]= \frac{1}{2}\beta\right\}=1,
\end{eqnarray*}
where last equality is a consequence of Lemma \ref{useful_prop}.

\hspace{12cm}$\Box$

\begin{cor}
\label{cor_3.1_no_frozen}
There exist $\tilde{N}>0$ and $M_0>0$ such that for every $N>\tilde{N}$
$$
\left\langle b^N(\xi,\Theta,x),\frac{x}{\vert x\vert}\right\rangle> \beta\gamma\left(\frac{N}{2}-1\right)\frac{1}{\vert x\vert} \quad \vert x\vert\geq M_0 \quad a.s.
$$
\end{cor}
\begin{flushleft}
Proof
\end{flushleft}
It is an immediate consequence of Proposition \ref{main_proposition_no_frozen}.

\hspace{12cm}$\Box$

\begin{prop}
\label{criterium_prop_no_frozen}
If $0<\beta\gamma\leq 1$, then equation (\ref{eq_z}) with $c>\frac{1}{\beta\gamma}\left(1+\frac{2}{\frac{N}{2}-1}\right)$ satisfies a.s. condition (\ref{cond_vere}) with $r=c\left(\frac{N}{2}-1\right)>\frac{N}{2}+1$.
\end{prop}
\begin{flushleft}
Proof
\end{flushleft}
It is an immediate consequence of Corollary \ref{cor_3.1_no_frozen} and Proposition \ref{prop_3.2}.

\hspace{12cm}$\Box$

In Proposition \ref{main_no_frozen} Veretennikov's Theorem is applied to equation (\ref{eq_z}), while the existence of an invariant measure for $y$ and the rate of convergence to it is derived in Corollary \ref{cor1_no_frozen}.

\begin{prop}
\label{main_no_frozen}
Let $\nu(t)$ be the distribution of $z(t)=(z_i(t))_{1\leq i\leq N}$ satisfying equation (\ref{eq_z}) with initial condition $z(0)$ and $c>\frac{1}{\beta\gamma}\left(1+\frac{2}{\frac{N}{2}-1}\right)$ with $0<\beta\gamma\leq 1$. For $N$ sufficiently large, there exists a.s. an invariant measure $\nu_{inv}$ for $z$ such that
$$
\vert \nu(t)-\nu_{inv}\vert\leq m(1+\vert z(0)\vert^l)(1+t)^{-(k+1)},
$$
where $m$ is a positive constant, $0<k<c\left(\frac{N}{2}-1\right)-\frac{N}{2}-1$ and $2k+2<l<2c\left(\frac{N}{2}-1\right)-N$.
\end{prop}
\begin{flushleft}
Proof
\end{flushleft}
It is a consequence of Proposition \ref{criterium_prop_no_frozen} and Theorem \ref{teo_vere}.

\hspace{12cm}$\Box$

\begin{cor}
\label{cor1_no_frozen}
Let $\nu'(t)$ be the distribution of $y(t)=(y_i(t))_{1\leq i\leq N}$ satisfying equation (\ref{eq_y}) with initial condition $y(0)$ and let $c=\frac{1}{\beta\gamma}\left(1+\frac{2}{\frac{N}{2}-1}+\frac{2}{\left(\frac{N}{2}-1\right)^2}\right)$ with $0<\beta\gamma\leq 1$. For $N$ sufficiently large, there exists a.s. an invariant measure $\nu'_{inv}$ for $y$ such that
\begin{equation}
\label{inequality_no_frozen}
\vert \nu'(t)-\nu'_{inv}\vert\leq m'(1+\vert cy(0)\vert^l)(1+t)^{-(k+1)},
\end{equation}
where $m'$ is a positive constant, $0<k<\left(\frac{N}{2}+1\right)\left(\frac{1}{\beta\gamma}-1\right)+\frac{1}{\beta\gamma}\frac{2}{\frac{N}{2}-1}$ and $2k+2<l<N\left(\frac{1}{\beta\gamma}-1\right)+\frac{2}{\beta\gamma}+\frac{2}{\beta\gamma}\frac{2}{\frac{N}{2}-1}$.
\end{cor}
\begin{flushleft}
Proof
\end{flushleft}
Since the limit as $N$ grows to infinity of $c$ is $\frac{1}{\beta\gamma}$ and $z=cy$, for $N$ sufficiently large $\vert \nu'(t)-\nu'_{inv}\vert=m'\vert \nu(t)-\nu_{inv}\vert$ and the thesis follows from Proposition \ref{main_no_frozen}.

\hspace{12cm}$\Box$

\subsection{Waiting time for reaching the stationary state}

Corollaries \ref{cor1} and \ref{cor1_no_frozen} provide a rate of convergence toward the invariant distribution for the score differences distribution; by making explicit $t$ from (\ref{inequality}) and (\ref{inequality_no_frozen}), it is hence possible to obtain the limiting behavior as $N$ grows to infinity of the waiting time for reaching the stationary state.
Being the criteria for the application of Veretennikov Theorem sufficient conditions for the existence of an invariant measure, the waiting time obtained may not be the smallest one, and hence we have to talk of upper bound of the waiting time for reaching the stationary state.

Since the scores time is the $y$'s time rescaled by $\Gamma$, i.e. $y_i(t)=\Gamma\frac{U_{+,i}(t/\Gamma)-U_{-,i}(t/\Gamma)}{2}$, in studying the waiting time for reaching a stationary state for the MG we have to refer to the time $\tau=\frac{t}{\Gamma}$, that is the own time of the scores $U_{s,i}$ corresponding to each agents possible choices $s=+1,-1$.

In Proposition \ref{waiting} the asymptotic behavior of the waiting time for reaching the stationary state both for a MG with maximally and finite asymmetric initial condition is obtained.

\begin{prop}
\label{waiting}
For every $\epsilon >0$ let $T$ be such that
$$
m'(1+\vert cy(0)\vert^l)(1+T\Gamma)^{-(k+1)}=\epsilon, 
$$
where $l$ and $k$ are as in Corollaries  \ref{cor1} and \ref{cor1_no_frozen} (for games with maximally and finite asymmetric initial condition respectively).

It follows that 
$$
\vert \nu'(T\Gamma)-\nu'_{inv}\vert\leq \epsilon
$$
and
\begin{itemize}
\item[i)]if $\vert y(0)\vert\to\infty$
\begin{equation}
\lim_{N\to\infty}\frac{T}{\vert y(0)\vert^2}=\frac{m'}{\epsilon\Gamma},
\end{equation}
\item[ii)]if $\lim_{N\to\infty}\frac{N'}{N}=\gamma$, where $0<\gamma\leq 1$, $N'=\#\left\{y_i : y_i(0)\neq 0\right\}$ and $\vert y(0)\vert<\infty$,
\begin{equation}
\label{waiting2}
\frac{m''}{\epsilon\Gamma}\leq\lim_{N\to\infty}\frac{T}{\vert y(0)\vert^2}\leq\lim_{N\to\infty}\frac{m''}{\epsilon\Gamma}\vert cy(0)\vert^{N\left(\frac{1}{\beta\gamma}-1\right)+\frac{2}{\beta\gamma}+\frac{2}{\beta\gamma}\frac{2}{\frac{N}{2}-1}-2},
\end{equation} 
\item[iii)]under same conditions of point {\it ii)} with $\beta\gamma=1$
\begin{equation}
\label{waiting3}
\lim_{N\to\infty}\frac{T}{\vert y(0)\vert^2}=\frac{m''}{\epsilon\Gamma},
\end{equation}
\end{itemize}
where $m'$ and $m''$ are positive constants.
\end{prop}
\begin{flushleft}
Proof
\end{flushleft} {\it i)}
By Corollary \ref{cor1}, for $N$ sufficiently large,
$$
T=\frac{\left(\frac{m'}{\epsilon}\left(1+\vert cy(0)\vert^l\right)\right)^{\frac{1}{k+1}}-1}{\Gamma},
$$
where $m'$ is a positive constant, $0<k<\frac{2}{\frac{N}{2}-1}$ and $2k+2<l<2+\frac{4}{\frac{N}{2}-2}$.

In the limit $N\to\infty$ we obtain $c=1$, $k=0$ and $l=2$; it follows that
\begin{equation}
\label{waiting1}
\lim_{N\to\infty}\frac{T}{\vert y(0)\vert^2}=\frac{m'}{\epsilon\Gamma}.
\end{equation}
{\it ii)} By Corollary \ref{cor1_no_frozen}
$$
\frac{\left(\frac{m''}{\epsilon}\left(1+\vert cy(0)\vert^{2k+2}\right)\right)^{\frac{1}{k+1}}-1}{\Gamma}<T<\frac{\left(\frac{m''}{\epsilon}\left(1+\vert cy(0)\vert^{N\left(\frac{1}{\beta\gamma}-1\right)+\frac{2}{\beta\gamma}+\frac{2}{\beta\gamma}\frac{2}{\frac{N}{2}-1}}\right)\right)-1}{\Gamma},
$$
and (\ref{waiting2}) follows.

{\it iii)} Since $0<k<\left(\frac{N}{2}+1\right)\left(\frac{1}{\beta\gamma}-1\right)+\frac{1}{\beta\gamma}\frac{2}{\frac{N}{2}-1}$ and $2k+2<l<N\left(\frac{1}{\beta\gamma}-1\right)+\frac{2}{\beta\gamma}+\frac{2}{\beta\gamma}\frac{2}{\frac{N}{2}-1}$,
in the limit $N\to\infty$ we obtain $c=1$, $k=0$, $l=2$ and (\ref{waiting3}) follows from (\ref{waiting2}).

\hspace{12cm}$\Box$

Since, if the initial condition is finite, $\vert y(0)\vert^2 = O(N)$ and since $\frac{P}{N}=\alpha>0$, from point {\it iii)} of Proposition \ref{waiting} it follows that $T=O(\frac{P}{\alpha\Gamma})$.



\section{Conclusions}
\label{conclusions}
By applying the Veretennikov's Theorem to the continuum time version of the MG, we have obtained an upper bound for the asymptotic behavior, as the number of agents grows to infinity, of the waiting time for reaching the stationary state in the asymmetric phase ($\alpha>\alpha_c$).

Since Veretennikov's Theorem gives a sufficient condition for the existence of an invariant measure and it applies only to stochastic differential equations with not null initial condition, the waiting time obtained may not be the smallest one (it is an upper bound) and it holds only for a MG with asymmetric initial condition.
The fewer are the agents with initial asymmetry in evaluating their strategies, the stronger must be their asymmetry: if their number is $o(N)$, then at least one agent must have maximally initial asymmetry ($\vert y_i(0)\vert=\infty$), while if their number is $O(N)$ the initial asymmetry must be inversely proportional to their fraction with respect to the agents population size ($\beta\gamma=1$). It follows that the single-agent's weakest initial asymmetry allowed by our result is $x_i$ such that $x_i\tanh(x_i)=1$ ($x_i\approx 1.2$), corresponding to a game where all the agents ($\gamma=1$) have asymmetric initial condition. Our result is simply not applicable to a game with initial asymmetry weaker than the previous one.
It is worth to note that the limit we have derived agrees with the rule of thumb in performing numerical simulations to wait a number of time steps proportional to $\frac{P}{\alpha\Gamma}$ in order to reach the stationary state and that, being $T=O(\frac{P}{\alpha\Gamma})$, the time to equilibrium increases as $\alpha>\alpha_c$ goes toward $\alpha_c$.

\section{Appendix}

Consider the $n$-dimensional stochastic differential equation
\begin{eqnarray*}
\label{vere}
dX(t)&=&b(X(t))dt+\sigma(X(t))dW(t)\\
X(0)&=&x \in \R^n,
\end{eqnarray*}
either with initial data $X(0)=x\in\R^n$.

$W(t)$ is an $m$-dimensional Wiener process with $m\geq n$, $b$ is a locally bounded Borel function from $\R^n$ with values on $\R^n$ and $\sigma$ a bounded continuous non-degenerate matrix $n\times m$-function and suppose that the drift term satisfies the following condition: there exist constants $M_0\geq 0$ and $r>0$ such that
\be
\label{cond_vere}
\left\langle b(x),\frac{x}{\vert x\vert}\right\rangle\leq -\frac{r}{\vert x\vert}, \quad \vert x\vert\geq M_0.
\ee

\begin{teo}[Veretennikov,\cite{veretennikov}]
\label{teo_vere}
Under assumption (\ref{cond_vere}) with $r>\frac{n}{2}+1$, for any $0<k<r-\frac{n}{2}-1$, with $l\in (2k+2,2r-n)$,
$$
\left\vert \mu_{x}(t)-\mu_{{\rm inv}}\right\vert\leq c(1+\vert x\vert^l)(1+t)^{-(k+1)},
$$
where $\left\vert \mu_{x}(t)-\mu_{{\rm inv}}\right\vert$ is the total variation distance between $\mu_{x}(t)$ and $\mu_{{\rm inv}}$, $c$ is a positive constant, $\mu_{x}(t)$ is the distribution of $X_t$, $x$ being the initial data, and $\mu_{{\rm inv}}$ is the invariant measure for $X_t$; in particular $\mu_{{\rm inv}}$ does exists.
\end{teo}



%
%
%
%

\end{document}